\documentclass[journal]{IEEEtran}
\IEEEoverridecommandlockouts   

\usepackage{amsthm}
\usepackage{amsmath,amssymb}
\usepackage{graphicx}
\usepackage{stmaryrd}
\usepackage{epstopdf}
\usepackage{xcolor}
\usepackage[english]{babel}
\definecolor{eff_blue}{RGB}{6,142,211}
\usepackage{tikz}
\usetikzlibrary{shapes}
\usepackage{adjustbox}
\usepackage[]{algorithm2e}
\usepackage{float}
\usepackage{eurosym}

\newcommand{\rpm}{\raisebox{.2ex}{$\scriptstyle\pm$}}

\makeatletter
\newcommand{\removelatexerror}{\let\@latex@error\@gobble}
\makeatother

\def\mathscr{\EuScript}

\newcommand{\cF}{\mathcal{F}}

\newcommand{\EE}{\mathbb{E}}

\newcommand{\RR}{\mathbb{R}}

\newcommand{\UU}{\mathbb{U}}

\newcommand{\WW}{\mathbb{W}}
\newcommand{\XX}{\mathbb{X}}

\newcommand{\np}[1]{(#1)}                                   % Parenth\`{e}se normal
\newcommand{\bp}[1]{\big(#1\big)}                           % Parenth\`{e}se big
\newcommand{\Bp}[1]{\Big(#1\Big)}                           % Parenth\`{e}se Big
\newcommand{\Bc}[1]{\Big[#1\Big]}                           % Crochet Big
\newcommand{\argmin}{\mathop{\arg\min}}                     % Arg-min
\newcommand{\espe}{\mathbb{E}}                              % Symbole esp\'{e}rance
\makeatletter
\def\va@a{\boldsymbol{\va@arg^{\textstyle\text{\unboldmath$\scriptstyle\va@expo$}}_{\textstyle\text{\unboldmath$\scriptstyle\va@index$}}}}
\def\va#1{\def\va@expo{}\def\va@index{}\def\va@arg{\uppercase{#1}}%
  \@ifnextchar^{\va@h}{\@ifnextchar_\va@u\va@a}}
\def\va@h^#1{\def\va@expo{#1}\@ifnextchar_\va@hu\va@a}
\def\va@u_#1{\def\va@index{#1}\@ifnextchar^\va@uh\va@a}
\def\va@hu_#1{\def\va@index{#1}\va@a}
\def\va@uh^#1{\def\va@expo{#1}\va@a}
\makeatother

\def\eqsepv{\; , \enspace}                                  % Virgule dans une \'{e}quation
\def\eqfinv{\; ,}                                           % Virgule en fin d'\'{e}quation
\def\eqfinp{\; .}                                           % Point en fin d'\'{e}quation
\newcommand{\finpreuvesymb}{$\Box$}%symbole fin preuve
\newcommand{\finremarksymb}{$\Diamond$}%symbole fin remarque
\newcommand{\finexemplesymb}{$\triangle$}%symbole fin exemple
\newcommand{\finpreuve}{\ \hspace*{\fill}\finpreuvesymb}
\newcommand{\finremark}{\ \hspace*{\fill}\finremarksymb}
\newcommand{\finexemple}{\ \hspace*{\fill}\finexemplesymb}
\makeatletter
\def\endproof{\finpreuve\@endtheorem}
\def\endremark{\finremark\@endtheorem}
\def\endexample{\finexemple\@endtheorem}
\makeatother

\begin{document}

\title{\LARGE Stochastic Optimization of Braking Energy Storage and Ventilation in a Subway Station}

\author{Tristan Rigaut, EFFICACITY, LISIS-COSYS-IFSTTAR, CERMICS-ENPC\\
Pierre Carpentier, UMA-ENSTA\\
Jean Philippe Chancelier, CERMICS-ENPC\\
Michel De Lara, CERMICS-ENPC\\
Julien Waeytens, LISIS-COSYS-IFSTTAR}

%\vspace{-0,4cm}
\maketitle
%\vspace{-0,4cm}

% \tableofcontents

\begin{abstract}
    %\vspace{-0.3cm}
In the Paris subway system, stations represent about one third of the overall energy
consumption. Within stations, ventilation is among the top consuming devices; 
it is operated at maximum airflow all day long, for air quality reasons. 
In this paper, we present a concept of energy system that displays comparable air
quality while consuming much less energy. The system comprises a battery that makes
it possible to recover the trains braking energy, arriving under the form of
erratic and strong peaks. We propose an energy management system (EMS) that, at
short time scale, controls energy flows and ventilation airflow. 
By using proper optimization algorithms, we manage to match supply with demand,
while minimizing energy daily costs. 
For this purpose, we have designed
algorithms that take into account the braking variability. They are based on the
so-called Stochastic Dynamic Programming (SDP) mathematical framework. 
We fairly compare SDP based algorithms with the widespread Model Predictive Control
(MPC) ones. 
First, both SDP and MPC yield energy/money operating savings of the order of one third,
compared to the current management without battery (our figure does not include
the cost of the battery).
Second, depending on the specific design, we observe that SDP outperforms MPC 
by a few percent, with an easier online numerical implementation.
\end{abstract}%

%\vspace{-0,6cm}

\section{Introduction}
\subsection{Context}
Apart from train traction, subway stations themselves represent a significant part (one third) of the energy
consumption of a subway system in cities like Paris. Fortunately, 
there is room to reduce their consumption by harvesting some of their unexploited energy
potential. We study here the potential energy recovery of a subway station
equipped with a battery to recover regenerative braking energy of subways. 

Ventilations are among the most significant energy consuming devices in subway
stations. One of the reason is because train braking produces a lot of particles
that need to be removed by ventilation. By producing regenerative braking
energy, trains can dissipate their kinetic energy with a lower bake pads
wear. Hence recovering braking energy improves air quality in stations. It might
then be useful to control simultaneously a ventilation and a battery to maximize
the benefits provided by this interaction. 

We present and compare hereby two classes of methods to solve optimal control problems in the presence of uncertainty. The first one is Model Predictive Control which requires only deterministic optimization tools. The second one is Stochastic Dynamic Programming based on Bellman equation. We apply two different flavours of SDP, a  state augmentation version to obtain the best possible performance we can achieve with SDP and a more classic one. We make a fair comparison of these methods by Monte Carlo simulation and present the results.

\subsection{Literature}
We survey literature on regenerative braking energy, air quality modelling
and energy storage management.

\subsubsection{Regenerative braking energy}
in most recent subway systems, trains already produce regenerative energy when they brake and transmit it to accelerating trains on the same line. However when there is no accelerating train nearby it is not possible to ensure the electrical supply demand equilibrium and regenerative braking is impossible. A comprehensive study of all possible ways to recover that energy is presented in \cite{gonzalez2013sustainable}. The authors of \cite{gonzalez2014systems} conclude that wayside energy storage is relevant to reduce the energy consumption of subway stations. The Southeastern Pennsylvania Transportation Authority (SEPTA) successfully installed wayside batteries to recover braking energy as reported in \cite{gillespie2014energy}. SEPTA is about to generalize the project to multiple stations.

\subsubsection{Air quality in buildings and subway stations}
an ANSES report \cite{anses} about air quality of underground subway stations states that the concentration of particulate matter whose size is inferior to $10 \mu m$ (PM10) can be unhealthy for the workers and maybe users. This is mainly due to ferrous PM10 that are generated during braking of the trains as stated in \cite{walther2017modelling} and \cite{grangepollution}. Subway stations operators in Paris took measures to monitor the concentration of PM10 \cite{grangepollution} that are openly available online. Many studies used Computational Fluid Dynamics technics to model the dispersion of pollutants in subway stations to produce predictive models as did the authors of \cite{camelli}. These methods are computationally very expensive and could hardly be integrated in an optimization problem without using reduced basis methods \cite{prudhomme:hal-00798326} that are challenging to implement in dynamic environments such as subway stations. The methods presented in \cite{walther2017modelling} and  \cite{doi:10.3155/1047-3289.58.4.502} use zonal models to compute an estimation of the global indoor air quality. These models are much more computationally efficient but require many approximations. The authors of \cite{vaccarini2016model} used MPC to control the energy consumption of ventilations and the related climate in a subway station. They estimate that their strategy could save up to $30\%$ of energy while maintaining the same comfort levels but don't manage an electrical storage simultaneously.

\subsubsection{Energy storage management}
most of the litterature apply MPC or Two Stage Stochastic Programming techniques to short term operation optimization of energy storage with uncertain supply as observed in \cite{papavasiliou2017application}. Authors of \cite{parisio2014model} and \cite{pflaum2014comparison} present MPC strategies to manage energy in battery and building climate.
In \cite{heymann2015continuous}, \cite{heymann2016stochastic} and  \cite{haessig2013computing} the authors present SDP strategies to control batteries in microgrids. In \cite{wu2016stochastic} SDP is applied to smart home management with electricle vehicle battery management. Few papers \cite{riseth2017comparison}, \cite{7100937} seem to compare the performance of different stochastic optimal control strategies.

\section{Energy system model}
\label{Energy_system_model}

We consider the energy system sketched in Figure~\ref{fig:elecschema}.
We present the equations describing its physical evolution in continuous time
(denoted by~$t$):
energy storage, Kirchoff laws and air quality. 
This energy system model will be the basis to simulate different 
management strategies corresponding to different EMS.

\subsection{Energy storage model}

We use a classical simple model of the dynamics of the energy storage system,
with the following variables:
\begin{itemize}
\item 
$s(t)$ (\%), the state of charge of the battery at time~$t$;
\item 
$u^b(t)$ ($kW$), the charge ($ u^b(t)\geq 0 $) or discharge ($ u^b(t) \leq 0 $) power of the battery at time $t$;
indeed, we observe on Figure~\ref{fig:elecschema} that the battery can draw
power on the noational grid or provide power to the station. 
\end{itemize}

The dynamics of the state of charge is\footnote{%
 We recall that $(x)^ += \max(x,0)$ and $(x)^- = \min(0,x)$.}:
\begin{equation}
\frac{ds}{dt} = \rho_c (u^b(t))^+ + \frac{1}{\rho_d} (u^b(t))^- \eqfinv 
\label{eq:socdyn}
\end{equation}
with charge/discharge efficiencies $\rho_c$ and $\rho_d$.
This simple linear dynamical model is relevant as long as we can ensure, 
by proper
management, that the state~$s(t)$ of charge is kept between proper bounds
\( \underline{s} \leq s(t) \leq \overline{s} \)
% \begin{equation}
% s_{min} \leq s(t) \leq s_{max} \eqfinp
% \end{equation} 
(like $30\% $ and $90\% $ of the capacity), which also ensures 
a good ageing of the battery.

\subsection{Kirchoff laws } % Supply/demand balance}

On Figure~\ref{fig:elecschema}, we observe that all flows must be balanced
at the central node, by Kirchoff laws.
The balance equation writes 
\begin{equation}
d(t) + u^v(t) + u^b(t) = b(t) + u^r(t) \eqfinp 
\label{eq:Balance}
\end{equation}
We comment the different terms:
\begin{itemize}
\item 
the station consumes a purely exogenous power $d(t)$ ($kW$) on the grid at time $t$;
\item 
the ventilations of the station consume a power $u^v(t)$ ($kW$);
this energy is controllable and we assume that it can be switched 
between two modes corresponding to two distinct airflows;
\item 
the trains produce a recoverable power $b(t)$ ($kW$) on the line;
\item 
the difference $u^{r}(t)=d(t) + u^v(t) + u^b(t) - b(t) $ is 
either the trains braking power in excess ($u^{r}(t) \leq 0$) that will be
wasted (if there is not enough demand, trains brake mechanically instead of electrically),
or the power in shortage ($u^{r}(t) \geq 0$) that will be drawn on the national grid to satisfy the demand of the station and, possibly, to charge the battery.
\end{itemize}

\tikzstyle{block} = [draw, fill=blue!20, rectangle, 
    minimum height=3em, minimum width=6em]

\begin{figure}
\begin{center}
\begin{adjustbox}{max totalsize={0.4\textwidth}{\textheight},center}
\begin{tikzpicture}[auto, node distance=3cm]
    % We start by placing the blocks
    \node [block, fill=white, draw, name=station] {\includegraphics[width=.15\textwidth]{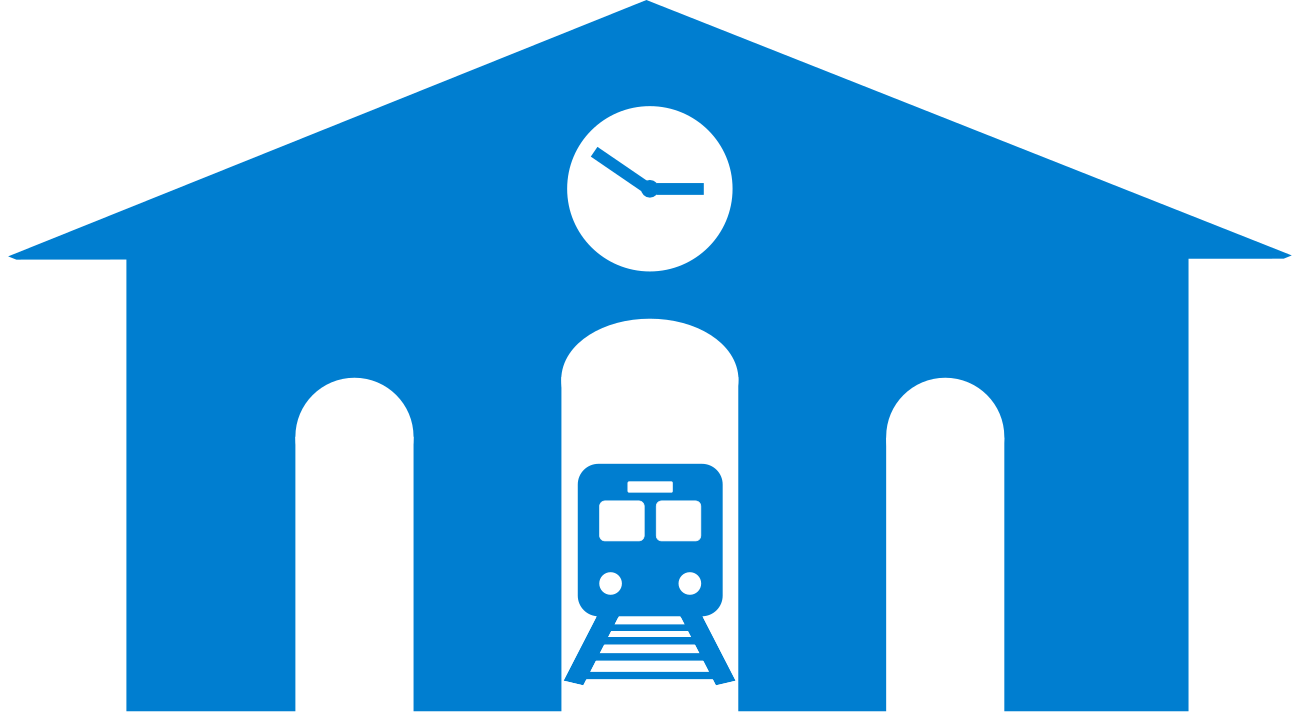}};
    \node [circle, draw, right of=station] (station_node) {};
    \node [label={$\va S$},block, fill=white, draw, right of=station_node] (battery) {\includegraphics[width=.12\textwidth]{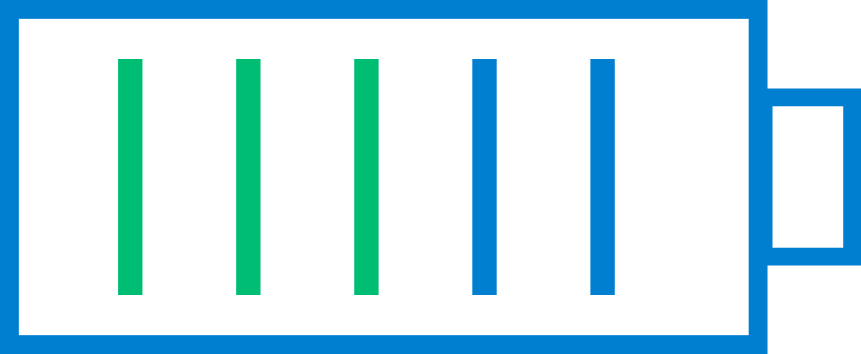}}; 
    \node [block, fill=white, draw, above of=station_node] (trains) {\includegraphics[width=.12\textwidth]{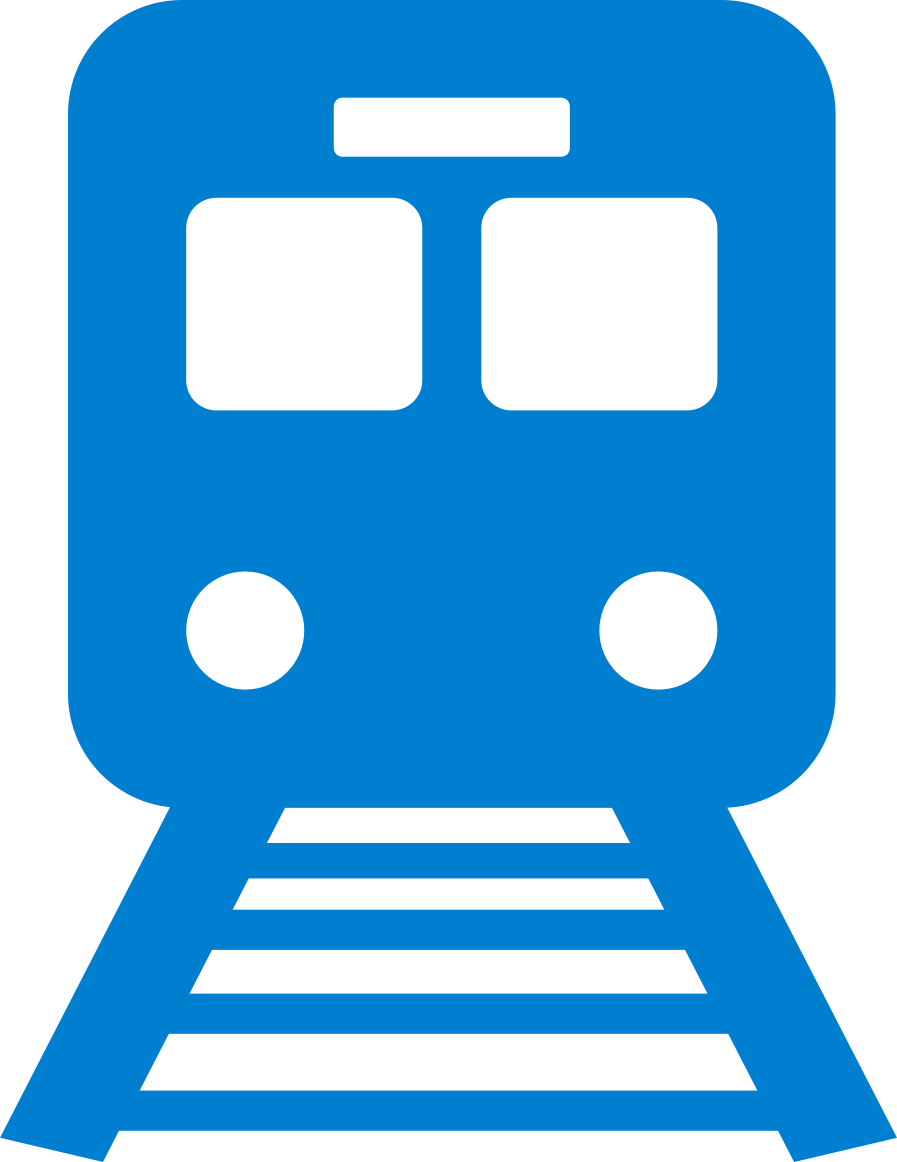}};
    \node [block, fill=white, draw, below of=station] (ventilation) {\includegraphics[width=.12\textwidth]{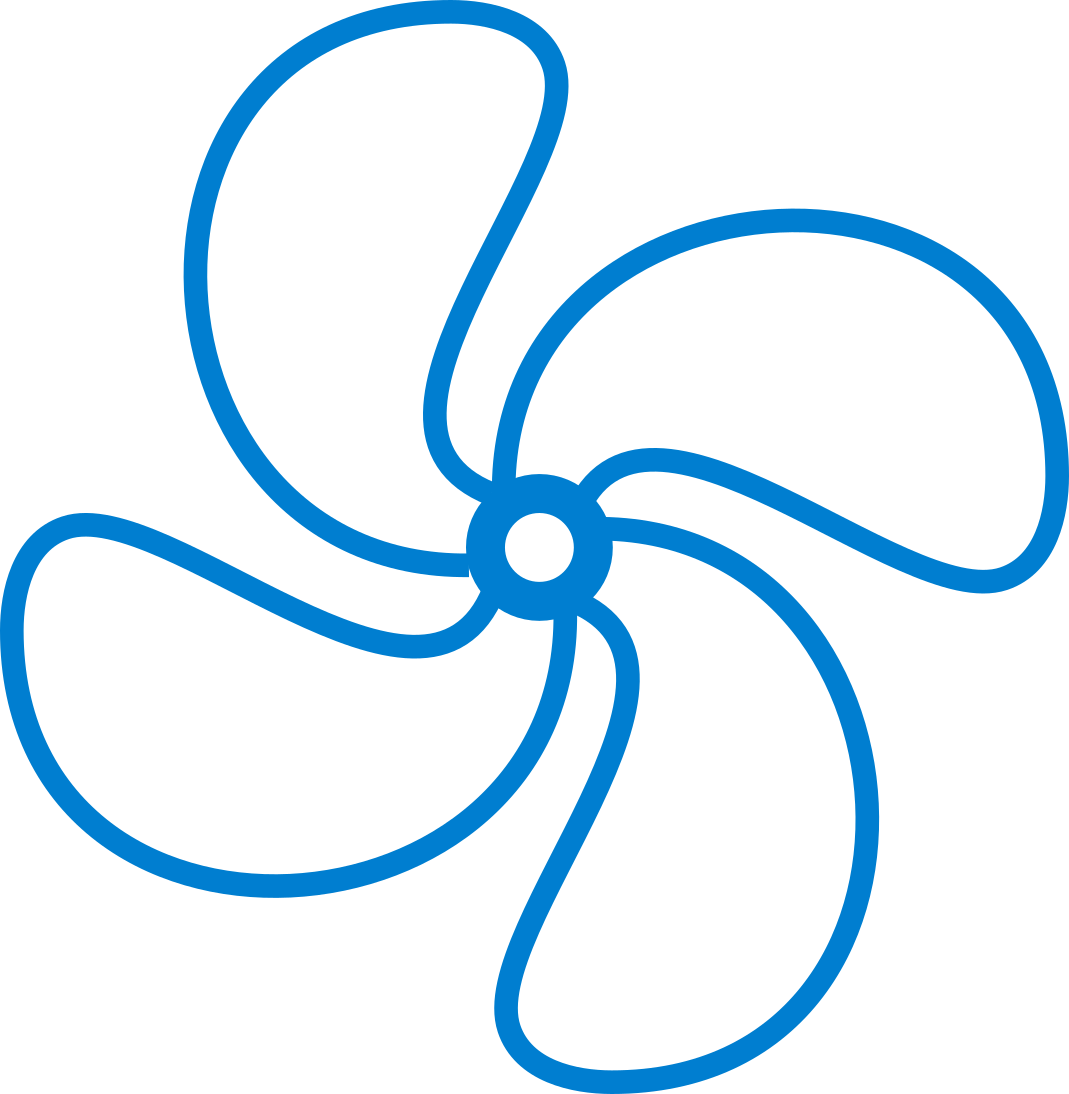}};
    
    \node [block, fill=white, draw, below of=battery] (national_grid) {\includegraphics[width=.1\textwidth]{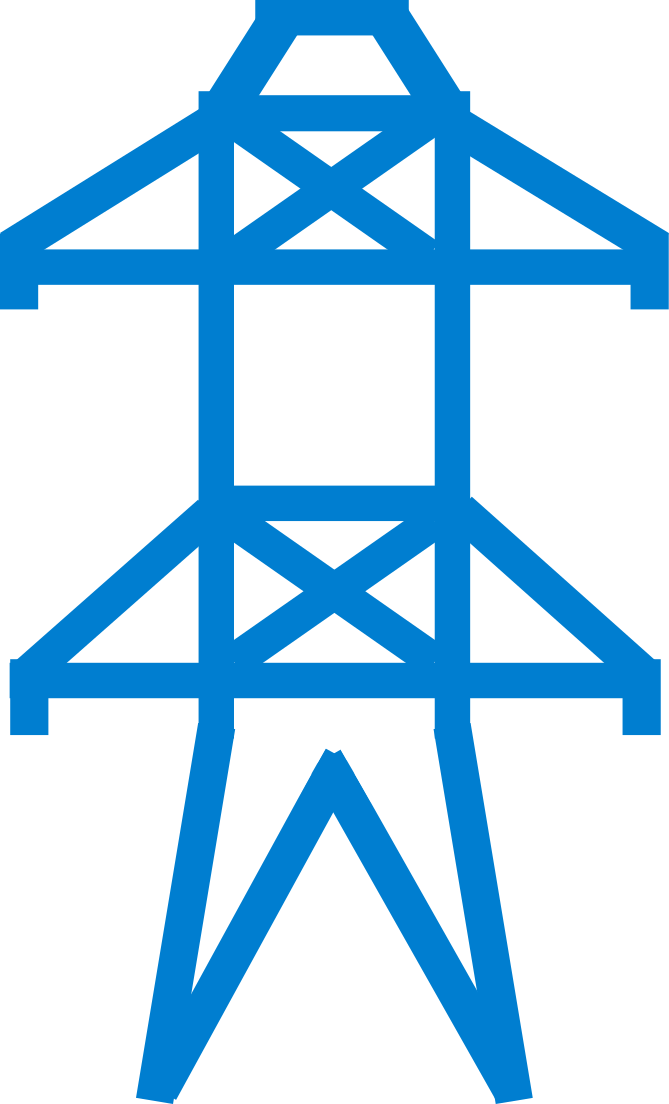}};

    \draw [draw,{<[scale=2]}-{>[scale=2]}] (station_node) -- node {$\va U^b$} (battery);
    \draw [draw,{<[scale=2]}-] (station) -- node {$\va D$} (station_node);
    \draw [draw,{<[scale=2]}-] (station_node) -- node {$\va B$} (trains);
    \draw [draw,{<[scale=2]}-] (station_node) -- node {$\va U^{r}$} (national_grid);
    \draw [draw,-{>[scale=2]}] (station_node) -- node {$\va U^v$} (ventilation);
\end{tikzpicture}
\end{adjustbox}
\caption{Electrical network representation}
\label{fig:elecschema}
\end{center}
\end{figure}
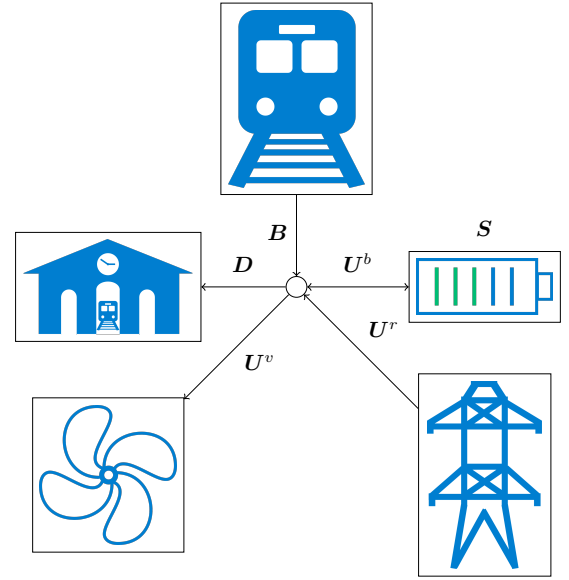

\subsection{Air quality model}

In \cite{doi:10.3155/1047-3289.58.4.502}, the authors use a bi-compartiment model in an office building in order to model the deposition/resuspension phenomenon. Due to the lack of data in subway stations to calibrate the surface dynamic model, we consider a simple air mass conservation model, as presented in \cite{walther2017modelling}, to model the dynamics of the particulate matters concentration in the subway station air. 
%In our case the lack of data to calibrate the surface dynamic model requires to make simplifications. 
As in \cite{walther2017modelling}, we assume that the floor is always saturated in dust particles so as to ignore the particles surface dynamics.
 Moreover, we assume that trains arriving in station produce particles 
by wearing brake pads and wheels, but also by resuspending particles from the floor. We use the model presented in \cite{walther2017modelling} to model the relation between trains arrivals and particles generation in the air. 

The dynamical equation for the PM10 concentration in the station is
\begin{align}
\frac{dc}{dt} &= \alpha n(t)^2 \nonumber  \\ &+ \Bp{ \frac{\rho_v}{v} u^v(t) + \beta n(t)} \Bp{c^o(t) - c(t)} -  \delta c(t) \eqfinv 
\label{eq:Cidyn}
\end{align}
with the following notations:
\begin{itemize}
\item 
$c(t)$ $(\mu g /m^3)$, PM10 concentration in the station air;
\item 
$c^o(t)$ $(\mu g /m^3)$, exogenous PM10 concentration outside the station;
\item 
$v$ $(m^3)$, volume of the station assimilated as a single zone;
\item 
$n(t)$ $(h^{-1})$, number of arriving trains per hour;
\item 
$\alpha$ $(\mu g h/m^3)$, apparent generation rate of particles by braking trains;
\item 
$\delta$ $(h^{-1})$, apparent deposition rate of particules;
\item 
$\beta$, apparent train contribution rate to natural ventilation;
\item 
$\rho_v$ $(m^3/kwh)$, global energy efficiency of the ventilations.
\end{itemize}

\subsection{Considerations on numerical simulations}
\label{numsim}
The equations~\eqref{eq:socdyn} and \eqref{eq:Cidyn}
form a system of ordinary differential equations.
We tested\footnote{%
We used the Julia \cite{bezanson2012julia} package DifferentialEquations.jl
  \cite{rackauckas2017differentialequations}.} 
that a forward Euler resolution with
$T_0 = 24h$ and $\Delta = 2~ \text{min}$ coincides with a 5th order Tsitouras
method using adaptative timestepping with a mean error of
 $0.06  \rpm 0.09 \%$. 
This makes it possible to simulate the energy system,
driven by given ventilation and battery control strategies, 
using a simple discrete time dynamical model.

\section{Optimization problem statement}

Once we dispose of the energy system dynamical model, 
we can envisage to simulate different 
management strategies and to compare them.
They are compared with respect to the daily costs that they induce, while
respecting constraints. To make this statement more formal and precise, 
we now formulate a mathematical optimization problem, under the form of a 
stochastic optimal control problem.

\subsection{Decisions are taken at discrete times}

By contrast with the energy system model developed
in Sect.~\ref{Energy_system_model}, where time is continuous, we adopt a
discrete time frame because decisions are made at discrete steps. 
Indeed, we consider a subway station grid equipped with a
hierarchical control architecture, as in most microgrids
\cite{olivares2014trends}, that needs time to compute and implement a decision. 
Decisions are produced every~$\Delta=2$~minutes, over an horizon $T_0 = 24h$;
then, they are sent to local controllers that make decisions at a faster pace. 

To make the connection with the variable indexed by continuous time
in Sect.~\ref{Energy_system_model}, we adopt the following convention:
for any variable~$x$, we put~$x_t = x(t\Delta)$ for $t = 0,\ldots,T =
\frac{T_0}{\Delta}$. In other words, $x_t$ denotes the value of the 
variable~$x$ at the beginning of the time interval \( [t, t+\Delta[ \). This dicretization is compatible with the one discussed in \S\ref{numsim}.

\subsection{Uncertainties are modelled as random variables}

% The electrical storage is used to mitigate production and demand uncertainty
% which means ensuring demand/production balance at every time without wasting too
% much the production or reducing the demand. 

We write random variables in capital bold letters, like $\va Z$, to distinguish them from deterministic variables $z$. 

We model energy demand~$\va D_t$ 
and trains braking energy production~$\va B_t$, defined 
when stating the balance equation~\eqref{eq:Balance}, as random variables. 
We do the same for the number~$\va N_t$ of trains
arrivals per hour and for the outside air quality~$\va C^o_t$,
both defined when stating the dynamical equation~\eqref{eq:Cidyn}
for the PM10 concentration in the station.

In the end, we define, for $t = 0,\ldots,T$, 
the vector of uncertainties at time step~$t$
\begin{equation}
\va W_t = (\va D_t, ~ \va B_t, ~ \va N_t, ~ \va C^o_t)^\top \eqfinp
\label{eq:W}
\end{equation}
We call $\va W_t$ the \emph{noise} at time~$t$, that is,
the uncertainties materialized at the end of the time interval~$[t-\Delta,t)$.
The noise $\va W_t$ takes value in the set  $\WW_t=\RR^{4}$.

\subsection{Control variables are modelled as random variables}

As time goes on, the noise variables $\va W_t$ are progressively unfolded 
and made available to the decision-maker. 
This is why, as decisions depend on observations in a stochastic optimal control problem, 
decision variables are random variables:
the variables in Sect.~\ref{Energy_system_model} will now become
random variables in capital bold letters.

At time step~$t$, at the beginning of the time interval \( [t, t+\Delta[ \),
the decision-maker takes two decisions: the battery charge/discharge power $\va U^b_t$
and the ventilation power $\va U^v_t$.
Then, at the end of the time interval \( [t, t+\Delta[ \)),
the decision-maker selects the power~$\va U^{r}_{t+1}$, drawn from the national grid, 
to react to the uncertainties~$\va D_{t+1}$ (demand) and $\va B_{t+1}$ (braking
energy) and to ensure the supply demand balance in the grid. 
This is made possible by a controlled DC/DC converter and supercapacitors 
that are not modelled in this problem. 
From the (balance equation) constraint \eqref{eq:Balance}:

\begin{equation}
\va U^{r}_{t+1} = \va D_{t+1} + \va U^v_t + \va U^b_t - \va B_{t+1} 
\eqfinp 
\label{eq:Et}
\end{equation} 
% We can therefore replace $\va U^{r}_{t+1}$ by its expression \eqref{eq:Et} in the objective \eqref{eq:Obj}. 
We group the two decision/control variables in a vector:
\begin{align}
\va U_t = (\va U^b_t, ~ \va U^v_t) \eqfinp
\end{align} 
We call $\UU_t=\RR^2$ the set in which the controls take their values.

\subsection{Non-anticipativity constraints for control variables}

To express the fact that the decision-maker (here the EMS) 
cannot anticipate on the future realizations of the noise,
we introduce $\cF_t$, 
the sigma algebra generated by all the past noises up to time $t$:
\begin{equation}
\cF_t = \sigma(\va W_0,\ldots,\va W_t) \eqfinp
\end{equation}
The increasing sequence $(\cF_0, \ldots, \cF_T )$ is the natural filtration used to model the information flow of the problem. 
The algebraic non-anticipativity constraint 
\begin{equation} \label{eq:mes}
\sigma(\va U_t) \subset \cF_{t} 
\end{equation} 
expresses the fact that the decision can only be made knowing no more 
than the past uncertainties \cite[chap. 4]{sowgbook}. 

We say that the controls satisfying \eqref{eq:mes} are $\cF_t$-measurable. 
%(???in a Decision-Hazard setting \cite{sowgbook}) 
Throughout the paper, a random variable $\va Z_t$
indexed by~$t$ is, by convention, $\cF_t$-measurable, that is,
$ \sigma(\va Z_t) \subset \cF_t$.

\subsection{State and dynamics}

In the energy system model developed
in Sect.~\ref{Energy_system_model},
the equations~\eqref{eq:socdyn} and \eqref{eq:Cidyn}
form a system of ordinary differential equations.
This is why we introduce two state variables,
the state of charge $s_t$ and the PM10 concentration $c_t$,
making thus a two-dimensional \emph{state} variable
\begin{align}
\va X_t = (\va S_t ,~ \va C_t)^\top \eqfinp
\end{align} 
We call $\XX_t=\RR^2$ the state space where the state
takes its values.

By sampling the continuous time differential 
equations~\eqref{eq:socdyn} and \eqref{eq:Cidyn}
at discrete time steps, 
and by considering that the control variables are piecewise constant
between two steps, we can define a discrete time dynamics
$f_t: \XX_t \times \UU_t \times \WW_{t+1} \to \XX_{t+1}$.
It is such that 
\begin{equation}
\va X_{t+1} = f_t(\va X_t, \va U_t, \va W_{t+1}) = 
\begin{pmatrix}
f^s_t(\va X_t, \va U_t, \va W_{t+1}) \\ 
f^c_t(\va X_t, \va U_t, \va W_{t+1}) 
\end{pmatrix}
\label{eq:dynamics_discrete_time}
\end{equation}
where 
\begin{subequations}
\begin{multline} \label{eq:sdyndis}
f^s_t(x_t, u_t, w_{t+1}) = s_t \\ 
+ \Delta \Bp{\rho_c (u^b_t)^+ + \rho_d^{-1}(u^b_t)^-} 
\end{multline}

\begin{multline}
f^c_t(x_t,u_t, w_{t+1}) = 
c_t -  \Delta \delta c_t + \Delta \alpha n_{t+1}^2   
\label{eq:cdyndis} \\ 
+\Delta \Bp{ \frac{\rho_v}{v} u^v_t + \beta n_{t+1}} 
\Bp{c^o_{t+1} - c_t} \eqfinp
\end{multline}  
\end{subequations}

\subsection{Bound constraints on the state and control variables}

As stated when writing the dynamics of the state of charge
in~\eqref{eq:socdyn}, the state of charge has to be kept bounded 
\begin{equation}
 \underline{s} \leq \va S_{t} \leq \overline{s} \eqfinp 
\label{eq:state_constraints}
\end{equation}
%As the dynamics~\eqref{eq:sdyndis} is not random, 
%this constraint turns into control constraints.

%Regarding the PM10 concentration, it is positive and 
%cannot reach unrealistically large values??????????????????:
%\begin{equation}
% c_{min} \leq \va C_{t} \leq c_{max} \eqfinp  ??????????????? inutile
%\end{equation}

The ventilation airflow can swith between two values,
 leading to the binary constraint
 \begin{equation}
   \va U^v_t \in \{ \underline{u^v},~ \overline{u^v} \} \eqfinv
\label{eq:control_constraints_binary}
 \end{equation}
and the charge/discharge power is limited,
leading to the  box constraint
\begin{equation}
  \underline{u^b} \leq \va U^b_{t} \leq \overline{u^b} \eqfinp
\label{eq:control_constraints_box}
\end{equation}

The bound
constraints~\eqref{eq:state_constraints}--\eqref{eq:control_constraints_binary}--\eqref{eq:control_constraints_box}, 
on the state and control variables,
can be summed in the synthetic expression 
\begin{equation}
  \np{\va X _t, \va U_t} \in B_t \subset \XX_t \times \UU_t \eqfinp 
\end{equation}

\subsection{The objective is an expected daily cost}

We consider the following criterion to be minimized:
\begin{equation}
\espe~ \Bc{
  \sum_{t=0}^{T-1} p_{t+1} \bp{\va U^{r}_{t+1}}^+ + \lambda \va C_{t+1}} 
\eqfinp 
\label{eq:Obj}
\end{equation}
We now comment each term. 

The term $\espe$ stands for the mathematical expectation. 
%As there is uncertainty we look for control strategies minimizing the criterion in expectation. With this \emph{risk measure} 
By the law of large numbers, mimimizing the mathematical expectation of costs
 ensures that the system will perform at its best over many days.

Inside the expectation, the sum over time represents the cumulated costs.
Those are a mix of two terms.

First, at every time step~$t$, we pay the electricity consumed on the national grid between $t-\Delta$ and $t$. We call $p_t$ $(\text{\euro{}}/kW)$ the cost of electricity per~$kW$ between $t-\Delta$ and $t$, that 
we assume to be deterministic. Therefore we pay 
\( p_t \times (\va U^{r}_{t})^+ ~(\text{\euro{}}) \) at time~$t$.
% \begin{equation}
% p_t \times \va U^{r}_{t} ~(\text{\euro{}}) \eqfinp \label{eq:Eco}
% \end{equation}

Second, we give a price of discomfort relative to air quality. 
Ideally, we would like to keep $T^{-1}\espe \Bp{\sum_{t=1}^{T} \va C_{t}}$,
the expected mean of particles concentration over a day, 
bounded. Indeed, this is the indicator used by the World Health Organization 
for its PM concentration guidelines \cite{who2006air}. 
To handle this constraint, we fix a marginal 
price~$\lambda$ $(\text{\euro{}}~m^3/\mu g))$  of discomfort 
associated with this ideal constraint. 
We have fixed this parameter by trials and errors,
after solving the problem for different values of $\lambda$. 
The cost of discomfort is then 
\( \lambda \times \va C_{t} ~(\text{\euro{}}) \).
% \begin{equation}
% \lambda \times \va C_{t} ~(\text{\euro{}}) \eqfinp \label{eq:Qai}
% \end{equation}

Finally, from~\eqref{eq:Obj} and~\eqref{eq:Et},
we define the \emph{instantaneous cost}
$L_t: \XX_t \times \UU_t \times \WW_{t+1} \to \RR$ by 
\begin{equation}
L_t(x_t, u_t, w_{t+1}) = 
p_t (d_{t+1} + u^v_t + u^b_t - b_{t+1})^+ + \lambda c_{t+1} 
\eqfinp
\label{eq:instantaneous_cost}
\end{equation}

\subsection{Stochastic optimal control problem formulation}

The EMS problem writes as a general Stochastic Optimal Control (SOC) \cite{sowgbook} problem in a risk neutral (expectation) setting
\begin{subequations}
  \label{eq:socproblem}
\begin{align}
  \underset{\va X,\va U}{\min} ~ & \EE \; \Bc{\sum_{t=0}^{T-1} L_t(\va X_t,\va U_t,\va W_{t+1})+K(\va X_T)} \label{eq:objective}\\
                              \text{s.t}~&\va X_{t+1} = f_t(\va X_t,\va U_t,\va W_{t+1}) \label{eq:dynamic}\\
&  \np{\va X _t, \va U_t} \in B_t \\                                 
                            & \sigma(\va U_t) \subset \cF_{t} \label{eq:nonanticipativity}
\end{align}
\end{subequations}
where $K$ is as final cost function --- which is $0$ in our case, as we are
indifferent of the state of charge at the end of the day. 

%Two classes of methods are considered herein to solve a SOC problem \cite{bertsekas1995dynamic} \cite{sowgbook} \cite{powell2007approximate}.

\section{Computation of online control strategies} 
\label{Online_computation}

%\subsection{Strategies as solutions of SOC problems}

The non anticipativity constraint~\eqref{eq:nonanticipativity} 
can be equivalently replaced by its functional counterpart~\cite[chap. 3, p86]{sowgbook}:
\begin{align}
\exists \pi_t: \WW_0 \times \ldots \times \WW_{t} \to \UU_t,~ \va U_t =
  \pi_t(\va W_0, \ldots, \va W_{t}) \eqfinp 
\label{eq:history_policy}
\end{align}
The mapping~$\pi_t$ is called a \emph{strategy} (more precisely a noise
dependent strategy).

In this paper, we restrict the search to solutions among the class 
of (augmented) \emph{state strategies} of the form
\begin{equation}
 \pi_t: \XX_t \times \WW_{t} \to \UU_t,~ 
\va U_t =   \pi_t(\va X_t, \va W_{t}) \eqfinp 
\end{equation}
This is indeed a restriction, as the state $\va X_t$ is, by the iterated 
dynamics~\eqref{eq:dynamics_discrete_time}, a function of 
\( (\va W_0, \ldots, \va W_{t}) \).
%In general we can search for functions of the physical state of the system and a subset of the past noises. We choose to keep only the state and the last past noise in this paper as the stochatic process we consider follows an AR1 process. 

In practice, we are not interested in knowing $\pi_t(x_t, w_{t}) $ for all possible
values of $(x_t, w_{t})$; we just want to be able to compute, on the fly, 
the value $u_t=\pi_t(x_t, w_{t}) $ when, at time~$t$, the couple $(x_t, w_{t})$
materializes.
This is why, in section~\ref{Model_Predictive_Control} and 
section~\ref{Stochastic_Dynamic_Programming},
we present two methods for the online implementation of strategies.
Both methods compute $u_t=\pi_t(x_t, w_{t}) $ by solving, online, a optimization problem.

\subsection{Model Predictive Control (MPC)}
\label{Model_Predictive_Control}

MPC is often casted in the context of deterministic optimization as it requires only to solve deterministic problems. However it can be often used to solve stochastic optimization problems. At time step~$t$, the MPC algorithm takes as inputs 
the state~$x$ of the system and all 
the previous uncertainties realizations $w_0,\ldots,w_{t}$.
One way or another, it selects a ``forecast'' 
$(\tilde w_{t+1}, \ldots, \tilde w_{T})$ and then solves the following
deterministic (open loop) optimal control problem:
\begin{subequations}
  \label{eq:detproblem}
\begin{align}
\underset{(u_{t}, \ldots, u_{T-1})}{\min} ~ & {\sum_{s=t}^{T-1} L_s(x_s,u_s,\tilde
                                              w_{s+1})+K(x_T)} 
\label{eq:MPC_objective}\\
                              \text{s.t}~&x_{s+1} = f_s(x_s,u_s,\tilde w_{s+1}) \\
& \np{x_s,u_s} \in B_s, x_t = x
\end{align}
\end{subequations}
From the optimal controls $(u_{t}, \ldots, u_{T-1})$ thus obtained,
the MPC algorithm only keeps the first $(\tilde u_{t}, \ldots, \tilde u_{t+N_{mpc}})$ 
(we call $N_{mpc}$ the \emph{reoptimization step} of the MPC).
Then, at time $t+N_{mpc}$, the MPC algorithm produces new controls 
by solving problem~\eqref{eq:detproblem} starting at $t+N_{mpc}$ 
with an updated forecast. 

As it proves delicate to select a decent forecast for all the remaining 
time horizon (and as a bad forecast can lead to poor decisions), 
the online problem horizon~$T-1$ in~\eqref{eq:MPC_objective}
is often cut at $t+h_t$, with $h_t \geq N_{mpc}$. 
Thus, one obtains problem \eqref{eq:detproblem} where 
the objective~\eqref{eq:MPC_objective} is replaced by 
\( \sum_{s=t}^{t+h_t} L_s(x_s,u_s,\tilde w_{s+1}) \).

%Finally the EMS runs the following algorithm at time $t$
%\begin{algorithm}
% \KwData{$t$, $w_0, \ldots, w_{t}, x_t$, $h_t$}
% \KwResult{$u_t^\sharp,\ldots,u_{t+\Delta_M}^\sharp $}
% Compute a forecast $(\tilde w_{t+1}, \ldots, \tilde w_{t+1+h_t})$\;
% Solve \eqref{eq:detproblemH} \;
% Return $u_t^\sharp,\ldots,u_{t+\Delta_M}^\sharp $\;
% 
% \caption{MPC: Online phase}
%\end{algorithm}
%
%
\subsection{Stochastic Dynamic Programming (SDP) based algorithms}
\label{Stochastic_Dynamic_Programming}
A major difference of MPC with the SDP methods is that there is no offline computation phase. 

\subsubsection{The offline-online SDPO algorithm encompasses two phases}
%The offline-online $SDP^{O}$ algorithm encompasses two phases:
a backward functional recursion performed offline;
a forward online optimization by exhaustive search.

Offline, the SDPO algorithm computes a sequence of functions~$\tilde{V}_t$ 
by backward induction as follows:
\begin{subequations}
\label{eq:Bellman}
\begin{alignat}{2}
  & \tilde{V}_T(x)= K(x)\\
 & \tilde{V}_t(x) = \underset{u \in \UU_t}{\min} && ~ \int_{\WW_{t+1}} \nonumber 
 \Bc{L_t(x,u,w_{t+1}) + \\ 
& && \tilde{V}_{t+1}\Bp{f_t(x,u,w_{t+1})} } \mu_{t+1}^{of}(dw_{t+1}) \eqfinp 
\label{eq:Bellmanobj}
\end{alignat}
\end{subequations}
Here, each $\mu_{t+1}^{of}$ is an (offline) probability distribution 
on the set~$\WW_{t+1}$. 
The recursion is often performed by exhaustive search in 
discretized versions of the state and control spaces,
hence requiring interpolation of the functions~$V_t$.
Indeed, $x_{t+1}=f_t(x,u,w)$ is not guaranteed to fall on a gridpoint of 
the discretized version of~$\XX_{t+1}$.

Online, at time~$t$, the SDPO algorithm uses the functions~$V_t$
and solves 
(with possibly a refined discretization of the control space~$\UU_t$)
\begin{align}
\nonumber u_t \in & \argmin_{u \in \UU_t} \int_{\WW_{t+1}} \Bc{L_t(x,u,w_{t+1}) + \\ 
& \tilde{V}_{t+1}\Bp{f_t(x,u,w_{t+1})} } \mu^{on}_{t+1}(w_{t}, dw_{t+1})  \eqfinp 
\end{align}
Here, $\mu_{t+1}^{on}$ is an (online) conditional probability distribution 
on the set~$\WW_{t+1}$, knowing the previous uncertainty~$w_{t}$. We choose a conditional distribution depending here only on the last uncertainty realization because we use an order 1 autoregressive model in our numerical experiment.
%a conditional probability distribution of $\va W_{t}$ that might benefit from our knowledge of the previous uncertainty realization~$w_{t}$ to narrow down future randomness. 
As the online conditional probability distribution $\mu_{t+1}^{on}$ 
depends on past uncertainties, this method produces state and noise dependent decisions in real time. 

%It leads to the following algorithm.
%\begin{figure}[H] 
%\removelatexerror
%\begin{algorithm}[H] \label{sdpoff}
% \KwData{Constants, random variables model}
% \KwResult{$\tilde{V}_{t=0,...,T}$}
% initialization $\tilde{V}_T = K$\;
% \For{$t=T-1...0$}{
%   \For{$x \in \XX_t$}{
%       $\tilde{V}_t(x) = +\infty$\;
%       \For{$u \in \UU_t$}{
%           $\tilde{V}_t(x) = \min \Bp{\tilde{V}_t(x), \underset{w \in \WW_{t+1}}{\sum} \mu_{t+1}(w) (L_t(x,u,w) + \tilde{V}_{t+1}(f_t(x,u,w)))}$
%       }
%   }
% }
% \caption{SDP: Offline phase}
%\end{algorithm}
%\end{figure}

%\begin{figure}[H] \label{sdpon}
%\removelatexerror
%\begin{algorithm}[H]
% \KwData{$t$, $\tilde{V}_{t+1}$, $w_0, \ldots, w_{t}, x_t$}
% \KwResult{$u_t^\sharp$}
% Quantize $\mu^o_t(w_0, \ldots, w_{t}, .)$\;
% \For{$u \in \UU_t$}{
%           $\underset{w \in \WW_{t+1}}{\sum} \mu^o_t(w_0,\ldots,w_t,w) (L_t(x,u,w) + \tilde{V}_{t+1}(f_t(x,u,w)))$
%       }
%Return best $u$\;
% \caption{SDP: Online phase}
%\end{algorithm}
%\end{figure}

It is well known \cite{bertsekas1995dynamic} that the above 
offline-online SDPO algorithm produces an optimal solution of the SOC
problem~\eqref{eq:socproblem} when i)
the random variables $\va W_0, \ldots, \va W_T$ are stagewise independent,
ii) $\mu_t^{of}$ is the probability distribution of $\va W_t$,
iii) $\mu_{t}^{on} = \mu_t^{of}$ is the (unconditional) probability distribution 
 of~$\va W_t$.

% \begin{subequations}
%   \label{eq:Bellman}
% \begin{align}
% & V_T(x)= K(x)\\
% & V_t(x) = \underset{u \in \UU_t}{\min} ~ \EE \; \Bc{ L_t(x,u,\va W_{t}) +
%   V_{t+1}\Bp{f_t(x,u,\va W_{t})} } 
% \label{eq:Bellmanobj}
% \end{align}
% \end{subequations} 
% \begin{align}
% \pi^\sharp_t(x) \in \argmin_{u \in \UU_t} \EE \; \Bc{ L_t(x,u,\va W_{t+1}) + V_{t+1}\Bp{f_t(x,u,\va W_{t+1})}}
% \end{align}
As, in our energy system case, the uncertainties are very likely correlated 
between successive time steps, they cannot be modelled by stagewise 
independent noises. Consequently, the strategy provided by the 
offline-online SDPO algorithm is not guaranteed to be optimal.

\subsubsection{The offline-online SDPA algorithm follows
the offline-online SDPO structure, but with}
the state~$x$ replaced by the couple $(x,w)$;
the uncertainty~$w$ replaced by a new uncertainty~$z$.

The dynamics~$f_t(x,u,w)$ is also replaced by 
a dynamics~$f^{A}_t\Bp{(x,w),u,z}$ of the form
\begin{align}
& f^{A}_t: \bp{\XX_t \times \WW_t} \times \UU_t \times \mathbb{Z}_{t+1} 
\to \bp{\XX_{t+1} \times \WW_{t+1}} \text{ where } \nonumber \\
&f^{A}_t\Bp{(x,w),u,z} = \Bp{f_t(x,u,f^w(w, z)),~f^w(w, z)} \eqfinp 
\end{align} 

It is straightforward that the above 
offline-online SDPA algorithm produces an optimal solution of the SOC
problem~\eqref{eq:socproblem} when 
there exists a stochastic process $\va Z_0,\ldots,\va Z_T$ such that
i) the random variables $\va Z_0, \ldots, \va Z_T$ are stagewise independent,
ii) \( \va W_{t+1} = f^w(\va W_{t}, \va Z_{t+1}) \).

The limit of this state augmentation strategy is the well known
\emph{curse of dimensionality}. The complexity of SDP grows exponentially 
with the number of state variables. Here, we try to handle a memory lag 
of one time step; but handling dependency between noises over 
multiple time steps would be out of reach.

\section{Numerical results, assessment and discussion}

In section \ref{Online_computation}, we outlined three methods to compute online
strategies. Now, we detail how to simulate them on the energy system model
developed in section~\ref{Energy_system_model} and how to compare their expected daily costs.

\subsection{Common data feeding the algorithms}

\subsubsection{Reference case}
we consider a subway station 
i) where the ventilation is operated at constant airflow $60~m^3/s$
ii) which is not equipped with a battery 
iii) which does not recover regenerative braking. 
With this ventilation strategy, the mean PM10 concentration over a day is $108~ \mu g /m^3$,
while the maximum is $182~ \mu g /m^3$. The consumption of the station over a day is $2.160 MWh$ which costs~$161$ \euro.

By choosing this reference case, our aim is to measure the daily savings 
made possible by investing into a battery
and by adopting one of the three strategies outlined 
in section \ref{Online_computation}.
This is a partial analysis, as we do not consider the costs of 
investment.

\subsubsection{Braking energy scenarios for algorithms design}
as stated in \eqref{eq:W}, the problem presents four sources of uncertainty. 
However, we assume that, in~\eqref{eq:W},
the demand~$\va D_t$,
the number~$\va N_t$ of trains per hour,
and the outdoor particles concentration~$\va C^o_t$ 
are deterministic in our numerical experiment.
Indeed, most of the uncertainty comes from the trains energy recovery
and, moreover, we can have pretty accurate forecasts for the variables 
that we assume deterministic.

A \emph{scenario} is any possible realization of the noise process 
$ (\va W_0, \ldots, \va W_T)$ written $ (w_0, \ldots, w_T)$.
For the braking energy, we generated $5,000$ so-called 
\emph{optimization scenarios} by using a rule, provided in the link in appendix A, 
calibrated on realistic data. 
These $5,000$ optimization scenarios are the common input provided 
to all the optimization algorithms, so that they can be used to design
the features of each algorithm.

\subsection{Numerical implementation of the MPC algorithm}

\subsubsection{Forecast}
\label{Forecast}
knowing a realization~$w_t$ of the noise $\va W_t$, 
we need to compute a forecast $ (\tilde{w}_{t+1}, \ldots, \tilde{w}_T)$ of the
future uncertainties.
The forecast relies upon the following log-AR$(1)$ model\footnote{%
The log transform ensures that we produce non negative forecasts.}
\begin{equation} 
\label{eq:log_AR1}
\log \va W_{t+1} = a \log \va W_{t} + \va Z_{t+1} \eqsepv\forall t = 0,\ldots, T \eqfinv
\end{equation}
with independent residual random variables $(\va Z_t)_{t=1,\ldots,T}$.
The coefficient~$a$ and the distribution of the residuals
are identified using the $5,000$ optimization scenarios.
%\begin{equation} 
% \label{eq:AR1}
% \va W_{t+1} = \va W_{t}^a \exp(\va Z_{t+1}) \eqsepv\forall t = 0,\ldots, T \eqfinp
% \end{equation}

\subsubsection{Deterministic problem resolution}

MPC requires to solve the determistic problem~\eqref{eq:detproblem}. We present the resolution method, based on a MILP formulation presented in Appendix~\ref{detpb_appendix}.

\subsection{Numerical implementation of the SDP algorithms}

\subsubsection{The offline-online SDPO algorithm}
it requires as input the probability distributions~$\mu_t^{of}$,
used to compute the functions~$V_t$ offline,
and the conditional probability distributions~$\mu_t^{on}$,
used to compute the controls online.
\begin{itemize}
\item 
$\mu_t^{of}$: we fit discrete probability distributions at each time step by
quantizing, using k-means algorithm, the values taken by the
$5,000$ optimization scenarios at this very time step~$t$.
\item 
$\mu_t^{on}$: knowing the realization~$w_{t-1}$, 
we obtain the conditional probability distributions~$\mu_t^{on}$ 
by using the formula $w_t = w_{t-1}^a\exp(z_t)$ (see~\eqref{eq:log_AR1}). 
From the $5,000$ optimization scenarios, we obtain 
$5,000$ values of~$z_t$, 
hence $5,000$ values of~$w_t$ by $w_t = w_{t-1}^a\exp(z_t)$. 
\end{itemize}

\subsubsection{The offline-online SDPA algorithm}
in addition to what is needed for the above SDPO algorithm,
it requires as input the new dynamics $f^w(w, z))$
such that \( \va W_{t+1} = f^w(\va W_{t}, \va Z_{t+1}) \). 
This dynamic is deduced from Equation~\eqref{eq:log_AR1}.

\subsection{Out of sample assessment of strategies}
We have generated $10,000$ so-called
\emph{assessment scenarios}, to be used \emph{only} for the assessment phase.

\begin{figure}
\begin{center}
\includegraphics[width=0.45\textwidth]{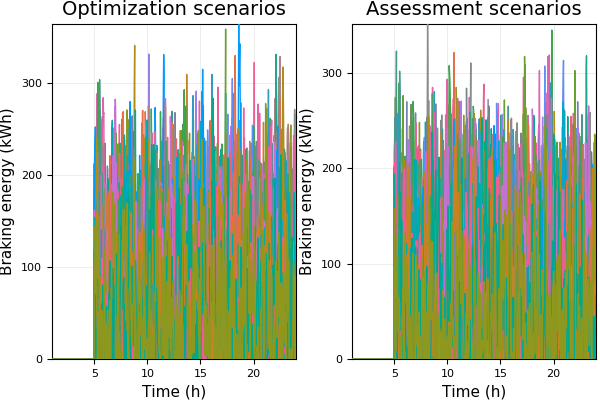}
\end{center}
\caption{Braking energy scenarios}
\label{brakingscens}
\end{figure}

We take good care to distinguish "optimization scenarios" from
"assessment scenarios", as displayed in Figure~\ref{brakingscens}. They are sealed.
Optimization scenarios were used to construct items entering the design
of the MPC and SDP algorithms. 
Assessment scenarios will be used to compare the strategies produced by 
these algorithms. This is what we call out of sample assessment. By this sealing, no algorithm can take advantage of the assessment scenarios
to be more fitted to the assessment phase.

The result of the assessment of a given strategy/algorithm
is an histogram of all the $10,000$ costs obtained 
along the assessment scenarios.

\subsection{Numerical results}

The computer used has 
Core i$7$, $2.2$Ghz processor
and $8$~Go ram + $12$~Go swap SSD memory.

\subsubsection{Comparing the algorithms performance}
the results\footnote{The lower the better, as we minimize costs. Results
  are $\pm$ the standard deviation.}  are summed up in Table~\ref{fig:restable}. We measure the savings with respect to a reference case with no regenerative braking recovered and ventilation at constant maximum speed over the day.

%\begin{figure}[H]
\begin{table}[H]
\begin{center}
\begin{tabular}{|l|c|c|c|l|} 
  \hline
  Strategy & SDPA & SDPO & MPC\\
  \hline
  ~ & ~ & ~ &\\
  Offline time & $3h47$ & $0h06$ & $0h00$ \\  
  Online time &  $4.0$~ms & $0.25$~ms & $8.5$~ms \\ 
  \hline
  ~ & ~ & ~ &\\
%  Objective & $316 \rpm 5.00$ & $319 \rpm 4.71$ & $321 \rpm 4.46$ \\ 
  Money savings (\euro) & $-74.1 \rpm 4.87$ & $-73.1 \rpm 4.54$ & $-71.1 \rpm 4.44$ \\
  PM10 $(\frac{\mu g}{ m^{3}})$ & $106 \rpm 0.11$ & $107 \rpm 0.11$ & $107 \rpm 0.08$\\
  Energy savings (kWh) & $-1050 \rpm 69.8$ & $-970 \rpm 59.55$ & $-942 \rpm 59.2$ \\
  \hline
\end{tabular}
\end{center}
\caption{Strategies performances comparison \label{fig:restable} }
\end{table}

We observe in Table~\ref{fig:restable} that all algorithms provide close
results. As we look in more detail, we see that SDPA outperforms both SDPO and
MPC on average for the economic savings, the mean PM10 concentration and the saved
energy. However, regarding the economic savings, the differences in mean performance 
(of order~3~\euro) are lower than standard deviations (of order~4.5~\euro), 
which makes it delicate to conclude.
The same analysis goes for the energy savings, although the confidence intervals
overlap less.

\begin{figure}[H]
\begin{center}
\includegraphics[width=0.45\textwidth]{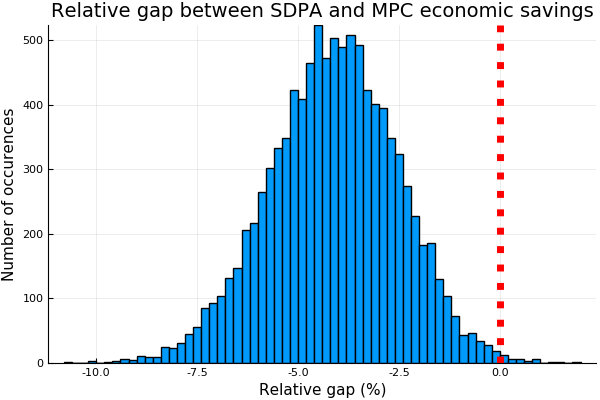}
\end{center}
\caption{Relative savings gap between SDPA and MPC}
\label{absgap}
\end{figure}

In fact, the three algorithms can be ranked as follows:
SDPA outperforms SDPO that outperforms MPC,
for the economic and energy savings (and they are comparable for air quality). 
To sustain this assertion, one has to look at Figure~\ref{absgap} 
that represents the distribution of the relative performance gap between MPC and
SDPA for the economic savings (a comparable analysis holds for the energy savings). 

%\begin{figure}[H]
%\begin{center}
%\includegraphics[width=0.45\textwidth]{relgap.png}
%\end{center}
%\caption{Relative objective gap between $SDP^{A}$ and $MPC$}
%\label{relgap}
%\end{figure}

On Figure~\ref{absgap}, the negative portion of the distribution to the left of the dotted red line is a
testimony in favor of SDPA. Our analysis of the assessment scenarios leads to
the following observations:
i) SDPA outperforms MPC for $9,967$ out of the $10,000$ scenarios, 
ii) SDPA outperforms SDPO for $8,221$ out of the $10,000$ scenarios, 
iii) SDPO outperforms MPC for all the scenarios. 

Concerning the computation time, Table~\ref{fig:restable} shows that SDPA
requires higher offline computation time than SDPO and MPC. 
As the online computation time for the three methods is way under $2$~minutes,
the three methods are implementable in real time (recall that the decision time
step is $2$~minutes). 
However, MPC differs from SDP algorithms along the following line:
MPC requires to solve a MILP online, so that there is no guarantee to reach the
optimum, or a feasible solution, within the prescribed $2$~minutes;
by contrast, both SDP algorithms only perfom an exhaustive search over all controls
in few milliseconds, which we consider safer for critical applications.

\subsubsection{Energy and air quality results}
we display and comment some energy and air quality results based on some of the $10,000$ assessment simulations.

\begin{figure}[H]
\begin{center}
\includegraphics[width=0.45\textwidth]{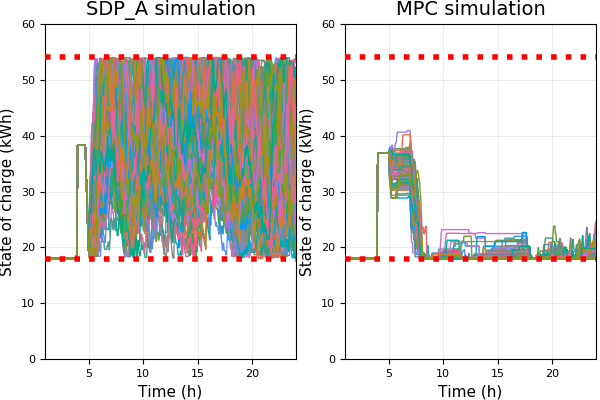}
\end{center}
\caption{Simulations of the state of charge}
\label{batfig}
\end{figure}

Figure~\ref{batfig} displays the state of charge trajectories of the battery on
the $10,000$ assessment scenarios for SDPA and MPC. We observe that the battery
is more intensively operated when using SDPA,
illustrating SDPA's ability to recover more energy than MPC.

\begin{figure}[H]
\begin{center}
\includegraphics[width=0.45\textwidth]{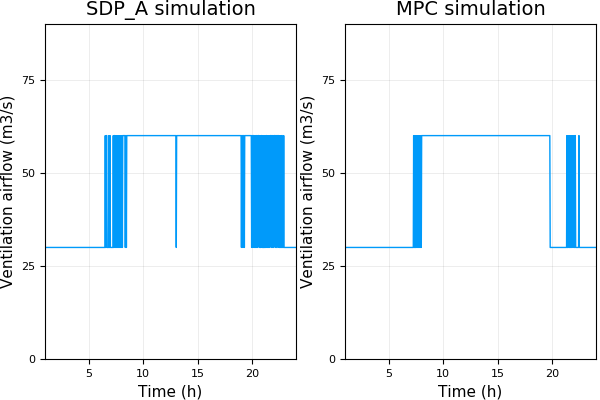}
\end{center}
\caption{Simulations of the ventilation airflow}
\label{ventilfig}
\end{figure}

\begin{figure}[H]
\begin{center}
\includegraphics[width=0.45\textwidth]{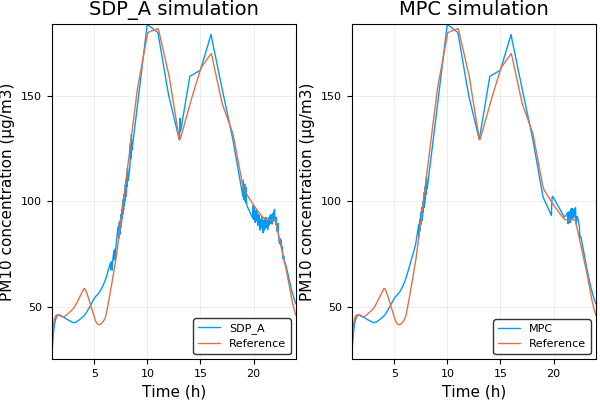}
\end{center}
\caption{Air quality simulations}
\label{qaifig}
\end{figure}

Figures \ref{ventilfig} and \ref{qaifig} display respectively the controls of the ventilation and the impact on the PM10 concentration over $1$~scenario for both SDPA and MPC. We recall that, in the reference case, the ventilation is operated at $60 m^3/s$ over the whole day. 
We observe that both algorithms decrease the consumption of the ventilation
while maintaining a similar air quality. \footnote{%
Had we modeled the particles generation reduction due to braking energy
recovery, we would have obtained a sharper decrease in~PM10 concentration.}

%\subsubsection{Discussion}
%we draw ??????????three conclusions from the numerical experiments.
%
%\begin{itemize}
%\item 
%The strength of MPC is that it can handle much more state variables than SDP. It is however limited to problems that can be efficiently solved by mathematical programming tools. SDP could still solve the problem with a continuous ventilation airflow while MPC would require to use nonlinear programming algorithms.
%\item 
%For critical applications SDP methods can be more robust than MPC as it requires to compute an exhaustive search with bounded maximum runtime while MPC could rely on NLP, MILP, MINLP solvers with uncertain runtime and uncertain quality of the solution.
%% \item 
%% Performances of the SDP strategies could be improved with a better discretization of the control and state spaces. 
%\end{itemize}

\section{Conclusions and perspectives}

We have presented a subway station energy system, with a battery recovering trains braking and smart control of
the ventilations. We have investigated methods to develop and implement an Energy Management
System that is able to handle uncertainties related to energy generation.
We have discussed the pros and cons of two popular techniques: 
Stochastic Dynamic Programming (SDP) and Model Predictive Control (MPC). 
For such a system (with a reasonable number of state variables), we have concluded 
that SDP is the best choice, even if MPC is a decent alternative. This is not the case in this paper but we recall that MPC
could require computationally expensive mathematical programming techniques 
to solve online deterministic problems. 
% For both methods, we note the importance of using the information 
% at our disposal online to take good decisions, 
% not just state feedbacks but past noises feedbacks. 

Our numerical experiments provide encouraging results. 
It seems that it pays to optimize to improve the energy efficiency 
and air quality of subway stations. Indeed, 
as seen on Figures~\ref{ventilfig} and~\ref{qaifig},
the ventilations energy consumption 
can be decreased without deteriorating the air quality. 

Our contribution is a first step towards the analysis of new subway station
energy systems.
It needs to be completed by an economic analysis that includes the costs of
batteries and the practical installation of such systems.

\appendix

\subsection{Data and parameters} \label{data_appendix} All the data used for the article is available on the following website:
https://trigaut.github.io/VentilationArticle.html

\subsection{Deterministic problem resolution} 
\label{detpb_appendix} 

To solve the MPC deterministic problem, we use mathematical programming
techniques by minimizing over states and control variables.
The dynamical equations~\eqref{eq:cdyndis}--\eqref{eq:sdyndis} are then simply equality constraints between decision variables. 
As the ventilation airflow~$u^v_t$ can switch between two modes, 
one of the decision variables is binary, 
leading to a Mixed Integer Non Linear Program (MILP). 
We use two simple tricks to turn~\eqref{eq:detproblem} 
into a Mixed Integer Linear Program. The constraint~\eqref{eq:cdyndis} contains the (non linear) product term
$u^v_t \times c_t$. To replace this term, we introduce the continuous 
variables~$a_t$ and the linear constraints 
\( 0 \leq a_t \leq \overline{C} \times u^v_t \) and
\( c_t - (1-u^v_t)\overline{C} \leq a_t \leq c_t \), 
% \begin{subequations}
% \begin{align*}
% & 0 \leq a_t \leq C_{max} \times u^v_t \eqfinv \\
% & c_t - (1-u^v_t)C_{max} \leq a_t \leq c_t \eqfinv
% \end{align*}
% \end{subequations} 
which ensures that $a_t = u^v_t \times c_t$ at optimality.

The constraint~\eqref{eq:sdyndis} contains positive and negative parts 
of~$u^b_t$, introducing non linearities. 
To circumvent the problem, we introduce two decision variables, 
$u^{b+}_t$ and $u^{b-}_t$ ($u^{b+}_t = (u^{b}_t)^+$ and $u^{b-}_t =
(u^{b}_t)^-$),
together with the constraint $u^{b+}_t \times u^{b-}_t = 0$.
It appears that this latter constraint can be removed as 
it always satisfied at optimality.
Indeed, there is no interest to flow through the battery to reach the demand 
as the battery efficiency coefficients waste power. 

To solve this MILP, we use the Julia package JuMP \cite{dunning2017jump} 
with the commercial solver Gurobi \cite{gurobi2015gurobi}.

%\vspace{-0,3cm}
%\nocite{6912056}
\footnotesize
\bibliographystyle{unsrt}
%\bibliography{bibliographic_data}
%\bibliographystyle{ieeetr} 
\bibliography{biblio}

\end{document}